\documentclass{amsart}

\usepackage{fancyhdr}
\usepackage{lastpage}
\usepackage{stmaryrd,yhmath}
\usepackage{tikz}

\newcommand{\bdis}{\begin{displaymath}}
\newcommand{\edis}{\end{displaymath}}
\newcommand{\be}{\begin{equation}}
\newcommand{\ee}{\end{equation}}
\newcommand{\mbb}{\mathbb}
\newcommand{\mcal}{\mathcal}

\newcommand{\vp}{\varphi}

\newcommand{\vth}{\vartheta}

\newcommand{\zf}{\zeta\left(\frac{1}{2}+it\right)}

\theoremstyle{definition}

\theoremstyle{remark}
\newtheorem{remark}[]{Remark}

\newtheorem*{mydef11}{{\bf Theorem 1}}

\newtheorem*{mydef12}{{\bf Theorem 2}}

\newtheorem*{mydef13}{{\bf Theorem 3}}

\newtheorem*{mydef21}{{\bf Definition 1}}

\newtheorem*{mydef22}{{\bf Definition 2}}

\newtheorem*{mydef23}{{\bf Definition 3}}

\newtheorem*{mydef24}{{\bf Definition 4}}

\newtheorem*{mydef25}{{\bf Definition 5}}

\newtheorem*{mydef26}{{\bf Definition 6}}

\newtheorem*{mydef41}{{\bf Corollary 1}}

\newtheorem*{mydef42}{{\bf Corollary 2}}

\newtheorem*{mydef51}{{\bf Lemma 1}}

\newtheorem*{mydef52}{{\bf Lemma 2}}

\newtheorem*{mydef53}{{\bf Lemma 3}}

\newtheorem*{mydef54}{{\bf Lemma 4}}

\newtheorem*{mydef55}{{\bf Lemma 5}}

\newtheorem*{mydef56}{{\bf Lemma 6}}

\newtheorem*{mydef81}{{\bf Property 1}}

\newtheorem*{mydef82}{{\bf Property 2}}

\newtheorem*{mydef83}{{\bf Property 3}}

\numberwithin{equation}{section}

\begin{document}

\title{Jacob's ladders, interactions between $\zeta$-oscillating systems and $\zeta$-analogue of an elementary trigonometric identity}

\author{Jan Moser}

\address{Department of Mathematical Analysis and Numerical Mathematics, Comenius University, Mlynska Dolina M105, 842 48 Bratislava, SLOVAKIA}

\email{jan.mozer@fmph.uniba.sk}


\begin{abstract}
In our previous papers, we have introduced within the theory of the Riemann zeta function the following notions: Jacob's ladders, oscillating systems,
$\zeta$-factorization, metamorphoses, \dots In this paper we obtain $\zeta$-analogue of an elementary trigonometric identity and other interactions
between oscillating systems.
\end{abstract}
\maketitle

\section{Introduction and survey of notions we have introduced in the theory of the Riemann zeta-function}

In our previous papers \cite{1} -- \cite{7} we have introduced within the theory of the Riemann zeta-function following notions: Jacob's ladders (JL),
$\zeta$-oscillating systems (OS), factorization formula (FF), metamorphosis of the oscillating systems (M), $\mcal{Z}_{\zeta,Q^2}$-transformation (ZT).

In the present paper we introduce the notion of interactions between oscillating systems (IOS).

The diagram of hierarchy of the previously mentioned notions is as follows:

\bdis
\mbox{JL}\rightarrow \mbox{OS} \rightarrow \mbox{FF} \rightarrow \left\{ \begin{array}{l}  \mbox{M}  \\ \mbox{ZT}  \\ \mbox{IOS} \end{array} \right.
\edis
where IOS represents new level of our theory.

The main result obtained in this direction is the following set of the $\zeta$-analogue of the elementary trigonometric identity
$\cos^2t+\sin^2t=1$:
\bdis
\begin{split}
& \cos^2(\alpha_0^{2,2})\prod_{r=1}^{k_2}
\left|\frac{\zeta\left(\frac 12+i\alpha_r^{2,2}\right)}{\zeta\left(\frac 12+i\beta_r^{2}\right)}\right|^2+
\sin^2(\alpha_0^{1,1})\prod_{r=1}^{k_1}
\left|\frac{\zeta\left(\frac 12+i\alpha_r^{1,1}\right)}{\zeta\left(\frac 12+i\beta_r^{1}\right)}\right|^2\sim 1,\ L\to\infty.
\end{split}
\edis
Of course, we will obtain also other interactions between oscillating systems.

\subsection{}

Let us remind that the Jacob's ladders
\bdis
\vp_1(t)=\frac 12\vp(t)
\edis
we have introduced in \cite{1} (see also \cite{2}), where the function $\vp(t)$ is an arbitrary solution to the nonlinear integral equation
\bdis
\int_0^{\mu[x(T)]}Z^2(t)e^{-\frac{2}{x(T)}t}{\rm d}t=\int_0^TZ^2(t){\rm d}t,
\edis
where each admissible function $\mu(y)$ generates the solution
\bdis
y=\vp(T;\mu)=\vp(T),\quad \mu(y)>7y\ln y.
\edis
We call the function $\vp_1(T)$ as Jacob's ladder as an analogue of the Jacob's dream in Chumash, Bereishis, 28:12.

\begin{remark}
By making use of these Jacob's ladders we have shown (see \cite{1}) that the classical Hardy-Littlewood integral (1918)
\bdis
\int_0^T\left|\zf\right|^2{\rm d}t
\edis
has - in addition to previously known Hardy-Littlewood expression (and other similar ones) possessing an unbounded error term at $T\to\infty$ - the
following infinite set of almost exact representations
\bdis
\begin{split}
& \int_0^T\left|\zf\right|^2{\rm d}t=\vp_1(t)\ln\{\vp_1(t)\}+ \\
& + (c-\ln 2\pi)\vp_1(t)+c_0+\mcal{O}\left(\frac{\ln T}{T}\right),\ T\to\infty ,
\end{split}
\edis
where $c$ is the Euler's constant and $c_0$ is the constant from the Titchmarsh-Kober-Atkinson formula.
\end{remark}

\subsection{}

Next, let us remind that we have introduced (see \cite{7}, (1.1)) the notion of the $(\zeta,Q^2)$-oscillating system. At this place we give the
definition in the following form.

\begin{mydef21}
The general oscillating system
\bdis
[\zeta,Q^2,k],\ k=1,\dots,k_0,\ k_0\in\mbb{N}
\edis
(with $k_0$ being arbitrary and fixed) is defined as follows
\bdis
[\zeta,Q^2,k]=\prod_{r=1}^k \left|\frac{\zeta\left(\frac 12+ix_r\right)}{\zeta\left(\frac 12+iy_r\right)}\right|^2,
\edis
where
\bdis
\begin{split}
& 0<T_0<x_1<x_2<\dots<x_k, \\
& T_0<y_1<y_2<\dots<y_k,
\end{split}
\edis
and
\bdis
\begin{split}
& (x_1,\dots,x_k),\ (y_1,\dots,y_k)\in (T_0,+\infty)^k \\
& x_r,y_r\not=\gamma:\ \zeta\left(\frac 12+i\gamma\right)=0,\ r=1,\dots,k
\end{split}
\edis
for sufficiently big $T_0$.
\end{mydef21}

\subsection{}

First of all we define the following class of admissible functions.

\begin{mydef22}
The symbol
\be \label{1.1}
f(t)\in \tilde{C}[T,T+U ]
\ee
stands for the following
\be \label{1.2}
f(t)\in C[T,T+U] \ \wedge \ f(t)>0
\ee
for
\be \label{1.3}
\begin{split}
& T>T_0,\ U\in (0,U_0], \\
& U_0=o\left(\frac{T}{\ln T}\right),\ T\to\infty.
\end{split}
\ee
\end{mydef22}

\begin{remark}
What concerns the last condition, see our paper \cite{3}, (7.1), (7.2).
\end{remark}

Next, we have shown in \cite{5}, (4.3) -- (4.18), (comp. \cite{7}, (2.1) -- (2.7)) the following is true:
\begin{itemize}
\item[(A)] there is a vector-operator $\hat{H}$ acting on $\tilde{C}$
\bdis
f(t)\mapsto \hat{H}f(t)=(\alpha_0,\alpha_1,\dots,\alpha_k,\beta_1,\dots,\beta_k),\ k=1,\dots,k_0
\edis
where
\be \label{1.4}
\begin{split}
& \alpha_r=\alpha_r(T,U,k;f),\ r=0,1,\dots,k , \\
& \beta_r=\beta_r(T,U,k),\ r=1,\dots,k,
\end{split}
\ee
\item[(B)] there is the factorization formula
\be \label{1.5}
\begin{split}
& f(t)\longrightarrow \prod_{r=1}^k \left|\frac{\zeta\left(\frac 12+ix_r\right)}{\zeta\left(\frac 12+iy_r\right)}\right|^2 = \\
& = \left\{ 1+\mcal{O}\left(\frac{\ln\ln T}{\ln T}\right)\right\}\frac{H(T,U;f)}{f(\alpha_0)}\sim \\
& \sim \frac{H(T,U;f)}{f(\alpha_0)},\ T\to\infty,
\end{split}
\ee
where
\bdis
H(T,U;f)=\frac{F(T+U;f)-F(T;f)}{U},
\edis
\bdis
F'_t(t;f)=f(t).
\edis
\end{itemize}

Now we give

\begin{mydef23}
Let
\bdis
f(t)\in\tilde{C}[T,T+U].
\edis
Then the particular oscillating system
\bdis
[\zeta,Q^2,k;f],\ k=1,\dots,k_0
\edis
is defined as follows
\be \label{1.6}
[\zeta,Q^2,k;f]=\prod_{r=1}^k \left|\frac{\zeta\left(\frac 12+i\alpha_r\right)}{\zeta\left(\frac 12+i\beta_r\right)}\right|^2 .
\ee
\end{mydef23}

\begin{remark}
In this case we have (see (\ref{3.10}), $c<0.6$)
\bdis
\begin{split}
& \alpha_{r+1}-\alpha_r>0.4\times \frac{T}{\ln T},\ 0,1,\dots,k-1 , \\
& \beta_{r+1}-\beta_r>0.4\times \frac{T}{\ln T},\ r=1,\dots,k-1,
\end{split}
\edis
and consequently, the inequalities of Definition 1 are fulfilled.
\end{remark}

\subsection{}

We have used words \emph{oscillating systems} in our Definitions 2 and 3. The main reason for this is in the
spectral form of the Riemann-Siegel formula (see \cite{5}, (3.1) -- (3.8))
\bdis
\begin{split}
& Z(t)=\sum_{n\leq \tau(x_r)}\frac{2}{\sqrt{n}}
\cos\{ t\omega_n(x_r)+\psi(x_r)\}+R(x_r), \\
& \tau(x_r)=\sqrt{\frac{x_r}{2\pi}}, \\
& R(x_r)=\mcal{O}(x_r^{-1/4}), \\
& t\in [x_r,x_r+V],\ V\in (0,x_r^{1/4}), \\
& Z(t)=e^{i\vth(t)}\zf \ \Rightarrow \ |Z(t)|=\left|\zf\right|,
\end{split}
\edis
where the functions
\bdis
\frac{2}{\sqrt{n}}
\cos\{ t\omega_n(x_r)+\psi(x_r)\}
\edis
are the Riemann's oscillators with:
\begin{itemize}
\item[(a)] the amplitude
\bdis
\frac{2}{\sqrt{n}},
\edis
\item[(b)] the incoherent local phase constant
\bdis
\psi(x_r)=-\frac{x_r}{2}-\frac{\pi}{8},
\edis
\item[(c)] the nonsynchronized local times
\bdis
t\in [x_r,x_r+V],
\edis
\item[(d)] local spectrum of cyclic frequencies
\bdis
\{\omega_n(x_r)\}_{n\leq \tau(x_r)},\ \omega_n(x_r)=\ln\frac{\tau(x_r)}{n}.
\edis
Similar formulae take place for $y_r$.
\end{itemize}

\begin{remark}
We have, for example, (see Definition 1)
\be \label{1.7}
\begin{split}
& [\zeta,Q^2,k]=\prod_{r=1}^k\left|\frac{\zeta\left(\frac 12+ix_r\right)}{\zeta\left(\frac 12+iy_r\right)}\right|^2= \\
& =\prod_{r=1}^k
\left|
\frac
{\sum_{n\leq \tau(x_r)}\frac{2}{\sqrt{n}}\cos\{ x_r\omega_n(x_r)+\psi(x_r)\}+R(x_r)}
{\sum_{n\leq \tau(y_r)}\frac{2}{\sqrt{n}}\cos\{ y_r\omega_n(y_r)+\psi(y_r)\}+R(y_r)}
\right|, \\
& k=1,\dots,k_0.
\end{split}
\ee
We see that (\ref{1.7}) is quite a complicated function (we can take $k_0$ arbitrarily big).
\end{remark}

\begin{remark}
The Riemann-Siegel formula (see \cite{8}, p. 60)
\bdis
Z(t)=\sum_{n\leq\tau(t)}\frac{2}{\sqrt{n}}\cos\{ \vth(t)-t\ln n\}+\mcal{O}(t^{-1/4}),
\edis
where
\bdis
\begin{split}
& Z(t)=e^{i\vth(t)}\zf, \\
& \vth(t)=-\frac t2\ln\pi+\mbox{Im}\ln\Gamma\left(\frac 14+i\frac t2\right)
\end{split}
\edis
is the Riemann formula that has been restored (by Riemann's manuscripts) and published by
C.L. Siegel.
\end{remark}

\begin{remark}
Let us notice that by our opinion the Riemann-Siegel formula represents Riemann's fundamental contribution to the theory
of oscillations (independently from the analytic number theory). Namely, the Riemann's oscillations are fated to describe of profound
laws of our Universe.
\end{remark}

\section{New formulae: interactions of oscillating systems and, especially, $\zeta$-analogue of the elementary
trigonometric identity}

\subsection{}

Let
\be \label{2.1}
\begin{split}
& f_1(t)=\sin^2t,\ t\in[\pi L+\mu,\pi L+\mu+U], \\
& L\in \mbb{N},\ U\in (0,U_0],\ \mu\geq \mu_0>0, \\
& 2\mu+U_0\leq \frac{\pi}{2}-\epsilon,\ \epsilon>0,\ \pi L+\mu>T_0
\end{split}
\ee
(with $\epsilon, \mu_0$ being sufficiently small and fixed). Of course,
\bdis
f_1(t)\in \tilde{C}[\pi L+\mu,\pi L+\mu+U] .
\edis
Now, if we use our algorithm (comp. third part of this paper) to the function $f_1(t)$, then we
obtain the following factorization formula
\be \label{2.2}
\prod_{r=1}^k\left|\frac{\zeta\left(\frac 12+i\alpha^1_r\right)}{\zeta\left(\frac 12+i\beta_r\right)}\right|^2\sim
\left\{\frac 12-\frac 12\frac{\sin U}{U}\cos(2\mu+U)\right\}\frac{1}{\sin^2(\alpha_0^1)},\ L\to\infty ,
\ee
where
\bdis
\begin{split}
& \alpha_r^1=\alpha_r(U,\mu,L,k;\sin^2t),\ r=0,1,\dots,k, \\
& \beta_r=\beta_r(U,\mu,L,k),\ r=1,\dots,k , \\
& \pi L+\mu<\alpha_0^1<\pi L+\mu+U \ \Rightarrow \ \mu<\alpha_0^1-\pi L<\mu+U, \\
& k=1,\dots,k_0.
\end{split}
\edis

\subsection{}

Next, we consider the function
\be \label{2.3}
f_2(t)=\cos^2t,\ t\in[\pi L+\mu,\pi L+\mu+U]
\ee
and obtain, by the similar way, the following factorization formula
\be \label{2.4}
\prod_{r=1}^k\left|\frac{\zeta\left(\frac 12+i\alpha^2_r\right)}{\zeta\left(\frac 12+i\beta_r\right)}\right|^2\sim
\left\{\frac 12+\frac 12\frac{\sin U}{U}\cos(2\mu+U)\right\}\frac{1}{\cos^2(\alpha_0^2)},\ L\to\infty ,
\ee
where
\bdis
\begin{split}
& \alpha_r^2=\alpha_r(U,\mu,L,k;\cos^2t),\ r=0,1,\dots,k, \\
& \mu<\alpha_0^2-\pi L<\mu+U,\qquad k=1,\dots,k_0.
\end{split}
\edis

\subsection{}

Consequently, we have obtained the following two sets
\be \label{2.5}
\begin{split}
& \left\{ [\zeta,Q^2,k_1;\sin^2t]\right\}_{k_1=1}^{k_0}, \\
& \left\{ [\zeta,Q^2,k_2;\cos^2t]\right\}_{k_2=1}^{k_0}
\end{split}
\ee
of particular oscillating systems.

\begin{remark}
We shall use the shortened phrase \emph{oscillating systems} instead of \emph{particular oscillating systems} in
similar cases.
\end{remark}

Now the following sets of new formulas is generated by two sets (\ref{2.5})
\be \label{2.6}
\begin{split}
& \cos^2(\alpha_0^{2,2})\prod_{r=1}^{k_2}\left|\frac{\zeta\left(\frac 12+i\alpha^{2,2}_r\right)}{\zeta\left(\frac 12+i\beta_r^2\right)}\right|^2+
\sin^2(\alpha_0^{1,1})\prod_{r=1}^{k_1}\left|\frac{\zeta\left(\frac 12+i\alpha^{1,1}_r\right)}{\zeta\left(\frac 12+i\beta_r^1\right)}\right|^2\sim 1,\ L\to\infty , \\
& k=1,\dots,k_0 ,
\end{split}
\ee
or, for example,
\be \label{2.7}
\begin{split}
& \prod_{r=1}^{k_2}\left|\frac{\zeta\left(\frac 12+i\alpha^{2,2}_r\right)}{\zeta\left(\frac 12+i\beta_r^2\right)}\right|^2\sim \\
& \sim \frac{1}{\cos^2(\alpha_0^{2,2})}-\frac{\sin^2(\alpha_0^{1,1})}{\cos^2(\alpha_0^{2,2})}
\prod_{r=1}^{k_1}\left|\frac{\zeta\left(\frac 12+i\alpha^{1,1}_r\right)}{\zeta\left(\frac 12+i\beta_r^1\right)}\right|^2,\ L\to\infty
\end{split}
\ee
(and second similar formula), where
\bdis
\begin{split}
& \alpha_0^{1,1}=\alpha_0^1(U,\mu,L,k_1;\sin^2t),\ \dots \\
& \beta_r^1=\beta_r(U,\mu,L,k_1), \\
& \alpha_0^{2,2}=\alpha_0^2(U,\mu,L,k_2;\cos^2t),\ \dots \\
& \beta_r^2=\beta_r(U,\mu,L,k_2).
\end{split}
\edis

\begin{remark}
We call the formula (\ref{2.6}) as the $\zeta$-analogue of the elementary trigonometric formula
\bdis
\sin^2t+\cos^2t=1.
\edis
\end{remark}

\begin{remark}
Of course, we have within (\ref{2.6}) huge number
\bdis
(k_0)^2
\edis
of distinct formulas.
\end{remark}

\subsection{}

It follows from (\ref{2.6}), (\ref{2.7}) that the corresponding oscillating systems (see (\ref{2.5}))
\be \label{2.8}
[\zeta,Q^2,k_1;\sin^2t],\ [\zeta,Q^2,k_2;\cos^2t]
\ee
are functionally depending systems. That is, if all parameters
\bdis
U,\mu,L,k_1,k_2
\edis
are fixed, then the corresponding values of the oscillating systems are linearly connected (in the asymptotic sense).

\begin{mydef24}
We shall call the above mentioned functional dependence of the oscillating systems (\ref{2.8}) as interaction between them and we shall
denote this by the following diagram
\be \label{2.9}
\begin{split}
& [\zeta,Q^2,k_1;\sin^2t] \longleftrightarrow [\zeta,Q^2,k_2;\cos^2t] , \\
& 1\leq k_1,k_2\leq k_0.
\end{split}
\ee
\end{mydef24}

Consequently, this paper is devoted to study of interactions between oscillating systems of the second order (as (\ref{2.9})), and also
to the third order systems, that is if
\bdis
[\zeta,Q^2,k_l,f_l] \longrightarrow (k_l,f_l),\ l=1,2,3.
\edis
In this sense, we will study the following diagrams
\bdis
(k_1,f_1) \longleftrightarrow (k_2,f_2),\ (k_1,f_1) \longleftrightarrow (k_2,f_2)\longleftrightarrow (k_3,f_3)\longleftrightarrow (k_1,f_1)
\edis
for
\bdis
1\leq k_1,k_2,k_3\leq k_0 .
\edis

\section{Short survey of our algorithm for generating the factorization formulae}

\subsection{}

If
\bdis
f(t)\in \tilde{C}[T,T+U]
\edis
then the formula
\be \label{3.1}
\begin{split}
& \frac 1U\int_T^{T+U}f(t){\rm d}t=H(T,U;f)>0,\ U\in (0,U_0] , \\
& H(T,U;f)=\frac{F(T+U;f)-F(T;f)}{U}, \\
& F_t'(t;f)=f(t)
\end{split}
\ee
holds true.

\subsection{}

Next, let us remind the we have introduced the following new type of integral in the theory of the Riemann zeta-function (see \cite{2}, (9.5), \cite{3}, (7.1), (7.2)): if
\bdis
f(t)\in \tilde{C}[T,T+U]
\edis
then
\be \label{3.2}
\begin{split}
& \int_{\overset{k}{T}}^{\overset{k}{\wideparen{T+U}}}f[\vp_1^k(t)]\prod_{r=0}^{k-1}\tilde{Z}^2[\vp_1^r(t)]{\rm d}t=\\
& = \int_T^{T+U}f(t){\rm d}t \ \Rightarrow \\
& \Rightarrow \ \int_{\overset{k}{T}}^{\overset{k}{\wideparen{T+U}}}\prod_{r=0}^{k-1}\tilde{Z}^2[\vp_1^r(t)]{\rm d}t=
U;\ f(t)=1,
\end{split}
\ee
where (see \cite{2}, (9.1), (9.2))
\be \label{3.3}
\begin{split}
& \tilde{Z}^2(t)=\frac{\left|\zf\right|}{\omega(t)}, \\
& \omega(t)=\left\{ 1+\mcal{O}\left(\frac{\ln\ln t}{\ln t}\right)\right\}\ln t,
\end{split}
\ee
and (see \cite{3}, (6.1), (6.4))
\bdis
[\overset{k}{T},\overset{k}{\wideparen{T+U}}];\
[T,T+U]\prec [\overset{1}{T},\overset{1}{\wideparen{T+U}}]\prec \dots \prec [\overset{k}{T},\overset{k}{\wideparen{T+U}}]
\edis
is the $k$-th reverse iteration of the segment $[T,T+U]$ where (see \cite{3}, (6.1))
\bdis
\vp_1(\overset{k}{T})=\overset{k-1}{T},\ \overset{0}{T}=T.
\edis

\subsection{}

Next, we obtain from (\ref{3.2}) by the mean-value theorem and (\ref{3.1}) (comp. \cite{5}, (4.5) -- (4.7))
\be \label{3.4}
\begin{split}
& f(\alpha_0)\prod_{r=1}^k \tilde{Z}^2(\alpha_r)=\frac{U}{\overset{k}{\wideparen{T+U}}-\overset{k}{T}}H(T,U;f), \\
& \alpha_r=\vp_1^{k-r}(d),\ r=0,1,\dots,k, \\
d=d(T,U,k;f),\ d\in \{ d\},
\end{split}
\ee
and, by the similar way, (comp. \cite{5}, (4.16), (4.17))
\be \label{3.5}
\begin{split}
& \prod_{r=1}^k \tilde{Z}^2(\beta_r)=\frac{U}{\overset{k}{\wideparen{T+U}}-\overset{k}{T}}, \\
& \beta_r=\vp_1^{k-r}(e),\ r=1,\dots,k, \\
e=e(T,U,k;f),\ e\in \{ e\} ,
\end{split}
\ee
where
\bdis
\{ d\}, \{ e\}
\edis
are sets of the abscises of the corresponding mean-values (\ref{3.4}) and (\ref{3.5}).

\begin{remark}
If we choose
\begin{itemize}
\item[(a)] only one
\bdis
d\in \{ d\} ,
\edis
\item[(b)]
only one
\bdis
e\in \{ e\} ,
\edis
\end{itemize}
then we obtain (see (\ref{3.4}), (\ref{3.5})) the corresponding vector-valued function
\be \label{3.6}
(\alpha_0,\alpha_1,\dots,\alpha_k,\beta_1,\dots,\beta_k),\ k=1,\dots,k_0,
\ee
i.e. we have $k_0$ corresponding vector-valued functions. Of course, these sets are defined for every admissible and fixed
\bdis
T,U,k;f(t).
\edis
\end{remark}

Now, we give the following

\begin{mydef25}
By making use the mean-value theorem three times (see (\ref{3.1}), (\ref{3.4}) and (\ref{3.5})) we define the vector-operator as follows
\be \tag{3.6}
\begin{split}
& \forall\- f(t)\in \tilde{C}[T,T+U] \longrightarrow \hat{H}f(t)= \\
& =(\alpha_0,\alpha_1,\dots,\alpha_k,\beta_1,\dots,\beta_k),\ k=1,\dots,k_0 .
\end{split}
\ee
\end{mydef25}

\begin{remark}
We notice explicitly that the operator $\hat{H}$ represents new type of operator since its definition is based on the function
\bdis
\zf
\edis
as well as on the Jacob's ladder (see (\ref{3.4}) and (\ref{3.5})).
\end{remark}

\subsection{}

Further we obtain, eliminating the factor
\bdis
\frac{U}{\overset{k}{\wideparen{T+U}}-\overset{k}{T}}
\edis
from (\ref{3.4}) and (\ref{3.5}), the following formula

\be \label{3.7}
\prod_{r=1}^k\frac{\tilde{Z}^2(\alpha_r)}{\tilde{Z}^2(\beta_r)}=\frac{H(T,U;f)}{f(\alpha_0)},\ T\to\infty .
\ee

\begin{remark}
We call formula (\ref{3.7}) as \emph{exact factorization formula}.
\end{remark}

Next, we obtain from (\ref{3.7}) by (\ref{3.3}), (comp. \cite{5}, (4.11) -- (4.14), \cite{7}, (2.7)) the following formula
\be \label{3.8}
\begin{split}
& \prod_{r=1}^k\left|\frac{\zeta\left(\frac 12+i\alpha_r\right)}{\zeta\left(\frac 12+i\beta_r\right)}\right|^2=
\left\{ 1+\mcal{O}\left(\frac{\ln\ln T}{\ln T}\right)\right\}\frac{H(T,U;f)}{f(\alpha_0)}\sim \\
& \sim \frac{H(T,U;f)}{f(\alpha_0)},\ T\to\infty.
\end{split}
\ee

\begin{remark}
We call formula (\ref{3.8}) as \emph{asymptotic factorization formula}.
\end{remark}

\subsection{}

In this paragraph, we present some properties of the vector operator $\hat{H}$ (properties of the corresponding vector-valued
function in (\ref{3.6})).

\begin{mydef81}
Since (see (\ref{1.2}), (\ref{1.3}))
\bdis
f(t)\in \tilde{C}[T,T+U] \ \Rightarrow \ f(t)>0,
\edis
then (see (\ref{3.1}), (\ref{3.3}) -- (\ref{3.5}))
\be \label{3.9}
\alpha_r\not=\gamma,\ \beta_r\not=\gamma:\ \left.\zf \right|_{t=\gamma}=0,\ r=1,\dots,k.
\ee
\end{mydef81}

\begin{remark}
Consequently, the set
\bdis
\{ \gamma\},\ \gamma>T_0
\edis
is the exceptional one for the last $2k$ components
\bdis
(\alpha_1,\dots,\alpha_k,\beta_1,\dots,\beta_k)
\edis
of the vector-valued function (\ref{3.6}).
\end{remark}

\begin{mydef82}
Next, the following inclusions holds true (comp. \cite{7}, (6.3))
\bdis
\alpha_0\in (T,T+U),\ \alpha_r,\beta_r\in (\overset{r}{T},\overset{r}{\wideparen{T+U}}),\ r=1,\dots,k.
\edis
\end{mydef82}

\begin{mydef83}
\be \label{3.10}
\begin{split}
& \alpha_{r+1}-\alpha_r\sim (1-c)\pi(T),\ r=0,1,\dots,k-1, \\
& \beta_{r+1}-\beta_r\sim (1-c)\pi(T),\ r=1,\dots,k-1,
\end{split}
\ee
where
\bdis
\pi(T)\sim \frac{T}{\ln T},\ T\to\infty
\edis
is the prime-counting function and $c$ is the Euler's constant (comp. \cite{5}, (2.8)).
\end{mydef83}

\begin{remark}
Jacob's ladder $\vp_1(T)$ can be viewed by our formula (see \cite{1}, (6.2))
\bdis
T-\vp_1(T)\sim (1-c)\pi(T)
\edis
as an asymptotic complementary function to the function
\bdis
(1-c)\pi(T)
\edis
in the following sense
\bdis
\vp_1(T)+(1-c)\pi(T)\sim T,\ T\to\infty .
\edis
\end{remark}

\begin{remark}
The asymptotic behavior of the sequences
\be \label{3.11}
\{ \alpha_r\}_{r=0}^k,\ \{ \beta_r\}_{r=1}^k
\ee
is by (\ref{3.10}) as follows (see \cite{5}, Remark 2): if $T\to\infty$ then the points of every sequence (\ref{3.11}) recede unboundedly
each from other and all together recede to infinity. Hence, at $T\to\infty$ each sequence in (\ref{3.11}) behaves as one-dimensional
Friedmann-Hubble universe.
\end{remark}

\section{First sets of lemmas}

\subsection{}

The following lemma holds true.

\begin{mydef51}
The function
\be \label{4.1}
f_1(t)=\sin^2t\in \tilde{C}[\pi L+\mu,\pi L+\mu+U],\ L\in\mbb{N},\ \pi L>T_0
\ee
corresponds with the following factorization formula
\be \label{4.2}
\begin{split}
& \prod_{r=1}^k\left|\frac{\zeta\left(\frac 12+i\alpha_r^1\right)}{\zeta\left(\frac 12+i\beta_r\right)}\right|^2\sim \\
& \sim \left\{ \frac 12-\frac 12\frac{\sin U}{U}\cos(2\mu+U)\right\}\frac{1}{\sin^2(\alpha_0^1)},\ L\to\infty,
\end{split}
\ee
where
\be \label{4.3}
\begin{split}
& \alpha_r^1=\alpha_r(U,\mu,L,k;\sin^2t),\ r=0,1,\dots,k, \\
& \beta_r=\beta_r(U,\mu,L,k),\ r=1,\dots,k, \\
& \pi L+\mu<\alpha_0^1<\pi L+\mu+U \ \Rightarrow \ \mu<\alpha_0^1-\pi L<\mu+U, \\
& k=1,\dots,k_0.
\end{split}
\ee
\end{mydef51}

\begin{proof}
Since
\bdis
\int \sin^2t{\rm d}t=\frac t2-\frac 14\sin 2t +C ,
\edis
then
\bdis
\int_{\pi L+\mu}^{\pi L+\mu+U} \sin^2t{\rm d}t=\frac U2-\frac 12\sin U\cos(2\mu+U),
\edis
and, of course,
\be \label{4.4}
\begin{split}
& \frac 1U\int_{\pi L+\mu}^{\pi L+\mu+U} \sin^2t{\rm d}t= \\
& = \frac 12-\frac 12\frac{\sin U}{U}\cos(2\mu+U).
\end{split}
\ee
Consequently, we obtain (\ref{4.2}) from (\ref{4.4}) by our algorithm (see third part of this paper).
\end{proof}

\subsection{}

By the similar way we obtain the following

\begin{mydef52}
The function
\be \label{4.5}
f_2(t)=\cos^2t\in \tilde{C}[\pi L+\mu,\pi L+\mu+U],\ L\in\mbb{N},\ \pi L>T_0
\ee
corresponds with the following factorization formula
\be \label{4.6}
\begin{split}
& \prod_{r=1}^k\left|\frac{\zeta\left(\frac 12+i\alpha_r^2\right)}{\zeta\left(\frac 12+i\beta_r\right)}\right|^2\sim \\
& \sim \left\{ \frac 12+\frac 12\frac{\sin U}{U}\cos(2\mu+U)\right\}\frac{1}{\cos^2(\alpha_0^2)},\ L\to\infty,
\end{split}
\ee
where
\be \label{4.7}
\begin{split}
& \alpha_r^2=\alpha_r(U,\mu,L,k;\cos^2t),\ r=0,1,\dots,k, \\
& \dots , \\
& \mu<\alpha_0^2-\pi L<\mu+U, \\
& k=1,\dots,k_0.
\end{split}
\ee
\end{mydef52}

\subsection{}

Now we obtain the following

\begin{mydef53}
The function
\be \label{4.8}
f_3(t)=\frac{1}{\cos^2t}\in \tilde{C}[\pi L+\mu,\pi L+\mu+U],\ L\in\mbb{N},\ \pi L>T_0
\ee
corresponds with the following factorization formula
\be \label{4.9}
\begin{split}
& \prod_{r=1}^k\left|\frac{\zeta\left(\frac 12+i\alpha_r^3\right)}{\zeta\left(\frac 12+i\beta_r\right)}\right|^2\sim \\
& \sim \frac{\sin U}{U}\frac{\cos^2(\alpha_0^3)}{\cos(\mu+U)\cos\mu},\ L\to\infty,
\end{split}
\ee
where
\be \label{4.10}
\begin{split}
& \alpha_r^3=\alpha_r(U,\mu,L,k;1/\cos^2t),\ r=0,1,\dots,k, \\
& \dots , \\
& \mu<\alpha_0^3-\pi L<\mu+U, \\
& k=1,\dots,k_0.
\end{split}
\ee
\end{mydef53}

\begin{proof}
Since
\bdis
\int \frac{{\rm d}t}{\cos^2t}=\tan t + C ,
\edis
then
\bdis
\int_{\pi L+\mu}^{\pi L+\mu+U} \frac{{\rm d}t}{\cos^2t}=\tan(\mu+U)-\tan\mu=\frac{\sin U}{\cos(\mu+U)\cos\mu},
\edis
and, of course,
\be \label{4.11}
\frac 1U\int_{\pi L+\mu}^{\pi L+\mu+U} \frac{{\rm d}t}{\cos^2t}=\frac{\sin U}{U}\frac{1}{\cos(\mu+U)\cos\mu}.
\ee
Consequently, the formula (\ref{4.9}) follows from (\ref{4.11}) by our algorithm.
\end{proof}

\section{Metamorphosis as the first interpretation of the factorization formula}

\subsection{}

First of all, we have (see (\ref{3.8})) that

\be \label{5.1}
\begin{split}
& \forall\- f(t)\in\tilde{C}[T,T+U] \ \rightarrow \\
& \rightarrow \ \prod_{r=1}^k\left|\frac{\zeta\left(\frac 12+i\alpha_r\right)}{\zeta\left(\frac 12+i\beta_r\right)}\right|^2\sim \\
& \sim \frac{H(T,U;f)}{f(\alpha_0)},\ T\to\infty,\ k=1,\dots,k_0, \\
& \alpha_r=\alpha_r(T,U,k;f),\ r=0,1,\dots,k , \\
& \beta_r=\beta_r(T,U,k;f),\ r=1,\dots,k ,
\end{split}
\ee
i.e. a factorization formula corresponds to every admissible function $f(t)$.

\subsection{}

Next, let us remind the general oscillating system (see Definition 1)
\be \label{5.2}
[\zeta,Q^2,k;f]=
\prod_{r=1}^k\left|\frac{\zeta\left(\frac 12+ix_r\right)}{\zeta\left(\frac 12+iy_r\right)}\right|^2,\ k=1,\dots,k_0,
\ee
where
\bdis
\begin{split}
& T_0<x_1<x_2<\dots<x_k , \\
& T_0<y_1<y_2<\dots<y_k, \\
& (x_1,\dots,x_k), (y_1,\dots,y_k)\in (T_0,+\infty)^k.
\end{split}
\edis
In connection with (\ref{5.1}) and (\ref{5.2}) we give the following.

\begin{remark}
We have introduced:
\begin{itemize}
\item[(a)] the notion of metamorphoses of an oscillating multiform in our paper \cite{4},
\item[(b)] the notion of metamorphoses of a quotient of two oscillating multiforms in our paper \cite{5}.
\end{itemize}
\end{remark}

\subsection{}

The mechanism of metamorphoses is as follows. If we get, by random sampling, such points
\be \label{5.3}
(x_1,\dots,x_k), (y_1,\dots,y_k)
\ee
that
\be \label{5.4}
(x_1,\dots,x_k)=(\alpha_1,\dots,\alpha_k);\quad
(y_1,\dots,y_k)=(\beta_1,\dots,\beta_k).
\ee
Then - at the points (\ref{5.3}) obeying property (\ref{5.4}) - the general oscillating system (\ref{5.2}) changes its old form
(chrysalis) into the new form (see (\ref{5.1}))
\bdis
\sim \frac{H(T,U,;f)}{f(\alpha_0)}
\edis
(butterfly), and the last is controlled by the function $\alpha_0$.

\section{$\mcal{Z}_{\zeta,Q^2}$-transformation as the second interpretation of the factorization formula}

\subsection{}

In \cite{7} we have introduced notion of $\mcal{Z}_{\zeta,Q^2}$-transformation. Namely, following transformation corresponds to the
factorization formula (\ref{5.1}): if
\bdis
f(t)\in \tilde{C}[T,T+U] ,
\edis
then
\be \label{6.1}
\begin{split}
& \begin{pmatrix} f(t) \\ t\in [T,T+U] \end{pmatrix}
\xrightarrow{\mcal{Z}_{\zeta,Q^2}}
\begin{pmatrix} \frac{H(T,U;f)}{f(\alpha_0)} \\ U\in (0,U_0] \end{pmatrix},\ T\to\infty .
\end{split}
\ee

\subsection{}

We may assume that the interpretation of (\ref{6.1}) is especially effective for such the signals (or pulses) that appear in the
theory of communication. In this direction we have introduced in our work \cite{7} the following.

\begin{mydef26}
The $\mcal{Z}_{\zeta,Q^2}$-transformation we call as $\mcal{Z}_{\zeta,Q^2}$-device with its input and output.
\end{mydef26}

In our paper \cite{7} we have obtained the following main law for the class of power-signals (pulses): in this case
\bdis
t^\Delta \in \tilde{C}[T,T+U],\ \forall\- \Delta\in\mbb{R},\ U\in (0,U_0],
\edis
\bdis
U_0<1,\ \forall\- L > L_0(\Delta),\ L\in\mbb{N}
\edis
we have (see (\ref{6.1}))
\be \label{6.2}
\begin{pmatrix} t^\Delta \\ t\in [L,L+U] \end{pmatrix}
\xrightarrow{\mcal{Z}_{\zeta,Q^2}}
\begin{pmatrix} 1 \\ U\in (0,U_0] \end{pmatrix},\ L\to\infty,
\ee
where
\bdis
\frac{H(L,U;f)}{(\alpha_0)^\Delta}\sim 1.
\edis

\begin{remark}
Every admissible power-signal (pulse)
\bdis
\begin{pmatrix} t^\Delta \\ t\in [L,L+U] \end{pmatrix}
\edis
on the input of the $\mcal{Z}_{\zeta,Q^2}$-device is transformed into telegraphic signal (pulse) that is into the unit
rectangular signal (pulse). For example,
\be \label{6.3}
\begin{pmatrix} t^{1000} \\ t\in [L,L+U] \end{pmatrix}
\xrightarrow{\mcal{Z}_{\zeta,Q^2}}
\begin{pmatrix} 1 \\ U\in (0,U_0] \end{pmatrix},\ L\to\infty,
\ee
\be \label{6.4}
\begin{pmatrix} t^{-1000} \\ t\in [L,L+U] \end{pmatrix}
\xrightarrow{\mcal{Z}_{\zeta,Q^2}}
\begin{pmatrix} 1 \\ U\in (0,U_0] \end{pmatrix},\ L\to\infty.
\ee
\end{remark}

\begin{remark}
We see that the unbounded signal (pulse) in (\ref{6.3}) as well as the negligible signal (pulse) in (\ref{6.4}) are both
transformed by $\mcal{Z}_{\zeta,Q^2}$-device into the telegraphic signal (pulse).
\end{remark}

\section{The $\zeta$-analogue of the elementary trigonometric identity as the first example of the interactions between
oscillating systems}

\subsection{}

Since all the coefficient-functions in the formulae (\ref{4.2}), (\ref{4.6}) are bounded and $\not=0$ (in the conditions
of (\ref{4.1})), then by eliminating the member
\bdis
\frac{\sin U}{U}\cos(\mu+U),
\edis
we obtain the following.

\begin{mydef11}
If the assumptions in (\ref{4.1}) are fulfilled then
\be \label{7.1}
\begin{split}
& \cos^2(\alpha_0^{2,2})\prod_{r=1}^{k_2}
\left|\frac{\zeta\left(\frac 12+i\alpha_r^{2,2}\right)}{\zeta\left(\frac 12+i\beta_r^2\right)}\right|^2+
\sin^2(\alpha_0^{1,1})\prod_{r=1}^{k_1}
\left|\frac{\zeta\left(\frac 12+i\alpha_r^{1,1}\right)}{\zeta\left(\frac 12+i\beta_r^1\right)}\right|^2\sim 1 , \\
& L\to\infty,\ 1\leq k\leq k_1,k_2\leq k_0,
\end{split}
\ee
where
\be \label{7.2}
\begin{split}
& \alpha_0^{1,1}=\alpha_0^{1}(U,\mu,L,k_1;\sin^2t), \dots \\
& \beta_r^1=\beta_r(U,\mu,L,k_1), \\
& \alpha_0^{2,2}=\alpha_0^{2}(U,\mu,L,k_2;\cos^2t), \dots \\
& \beta_r^2=\beta_r(U,\mu,L,k_2) .
\end{split}
\ee
\end{mydef11}

\begin{remark}
We call the formula (\ref{7.1}) as the $\zeta$-analogue of the elementary trigonometric identity
\bdis
\cos^2t+\sin^2t=1.
\edis
Of course, the formula (\ref{7.1}) denotes general element of the set
\bdis
(k_0)^2
\edis
distinct formulas for every admissible and fixed $L$ and for every fixed segment
\bdis
[\pi L_\mu,\pi L+\mu+U].
\edis
\end{remark}

\begin{remark}
Further, the formula (\ref{7.1}) expresses the asymptotic dependence of every pair of oscillating systems
\bdis
[\zeta,Q^2,k_1;\sin^2t],\ [\zeta,Q^2,k_2,\cos^2t],\quad 1\leq k_1,k_2\leq k_0
\edis
(comp. third part of this paper).
\end{remark}

\subsection{}

Now, we give explicitly (comp. (\ref{2.7}) and Definition 4) the following (see (\ref{7.1})).

\begin{mydef41}
\be \label{7.3}
\begin{split}
& \prod_{r=1}^{k_2}\left|\frac{\zeta\left(\frac 12+i\alpha_r^{2,2}\right)}{\zeta\left(\frac 12+i\beta_r^2\right)}\right|^2\sim \\
& \sim \frac{1}{\cos^2(\alpha_0^{2,2})}-\frac{\sin^2(\alpha_0^{1,1})}{\cos^2(\alpha_0^{2,2})}
\prod_{r=1}^{k_1}\left|\frac{\zeta\left(\frac 12+i\alpha_r^{1,1}\right)}{\zeta\left(\frac 12+i\beta_r^1\right)}\right|^2,\ L\to\infty,
\end{split}
\ee
\be \label{7.4}
\begin{split}
& \prod_{r=1}^{k_1}\left|\frac{\zeta\left(\frac 12+i\alpha_r^{1,1}\right)}{\zeta\left(\frac 12+i\beta_r^1\right)}\right|^2\sim \\
& \sim \frac{1}{\sin^2(\alpha_0^{1,1})}-\frac{\cos^2(\alpha_0^{2,2})}{\sin^2(\alpha_0^{1,1})}
\prod_{r=1}^{k_2}\left|\frac{\zeta\left(\frac 12+i\alpha_r^{2,2}\right)}{\zeta\left(\frac 12+i\beta_r^2\right)}\right|^2,\ L\to\infty .
\end{split}
\ee
\end{mydef41}

Let us remind (comp. (\ref{2.5})) that in our case, we have two sets of oscillating systems
\be \label{7.5}
\begin{split}
& \left\{ [\zeta,Q^2,k_1;f_1]\right\}_{k_1=1}^{k_0},\ f_1=f_1(t)=\sin^2t, \\
& \left\{ [\zeta,Q^2,k_2;f_2]\right\}_{k_2=1}^{k_0},\ f_2=f_2(t)=\cos^2t .
\end{split}
\ee

\begin{remark}
By the formulas (\ref{7.3}) and (\ref{7.4}) is expressed the property that we call (see Definition 4) as the interaction
between two corresponding oscillating systems from distinct sets (\ref{7.5}). We use the following diagram
\be \label{7.6}
[\zeta,Q^2,k_1;f_1] \longleftrightarrow [\zeta,Q^2,k_2,f_2]
\ee
to denote mentioned interaction.
\end{remark}

\section{The second case: interactions between three oscillating systems}

\subsection{}

Since (see (\ref{4.2}), (\ref{4.6}))
\be \label{8.1}
\begin{split}
& \frac{\cos^2(\alpha_0^{2,2})}{\cos(2\mu+U)}
\prod_{r=1}^{k_2}\left|\frac{\zeta\left(\frac 12+i\alpha_r^{2,2}\right)}{\zeta\left(\frac 12+i\beta_r^2\right)}\right|^2-
\frac{\sin^2(\alpha_0^{1,1})}{\cos(2\mu+U)}
\prod_{r=1}^{k_1}\left|\frac{\zeta\left(\frac 12+i\alpha_r^{1,1}\right)}{\zeta\left(\frac 12+i\beta_r^1\right)}\right|^2\sim \\
& \sim \frac{\sin U}{U},\ L\to \infty
\end{split}
\ee
then we obtain (see (\ref{4.9})) the following

\begin{mydef12}
If the assumptions of (\ref{4.1}) are fulfilled then
\be \label{8.2}
\begin{split}
& \prod_{r=1}^{k_3}\left|\frac{\zeta\left(\frac 12+i\alpha_r^{3,3}\right)}{\zeta\left(\frac 12+i\beta_r^3\right)}\right|^2\sim \\
& \sim \frac{\cos^2(\alpha_0^{2,2})\cos^2(\alpha_0^{3,3})}{\cos(2\mu +U)\cos(\mu+U)\cos\mu}
\prod_{r=1}^{k_2}\left|\frac{\zeta\left(\frac 12+i\alpha_r^{2,2}\right)}{\zeta\left(\frac 12+i\beta_r^2\right)}\right|^2- \\
& -\frac{\sin^2(\alpha_0^{1,1})\cos^2(\alpha_0^{3,3})}{\cos(2\mu +U)\cos(\mu+U)\cos\mu}
\prod_{r=1}^{k_1}\left|\frac{\zeta\left(\frac 12+i\alpha_r^{1,1}\right)}{\zeta\left(\frac 12+i\beta_r^1\right)}\right|^2,\ L\to\infty ,
\end{split}
\ee
where (comp. (\ref{4.10}))
\be \label{8.3}
\begin{split}
& \alpha_0^{3,3}=\alpha_0^3(U,\mu,L,k_3,f_3), \dots \\
& \beta_r^3=\beta_r(U,\mu,L,k_3), \\
& f_3=f_3(t)=\frac{1}{\cos^2t}
\end{split}
\ee
and for the symbols
\bdis
\alpha_0^{2,2},\ \alpha_0^{1,1} , \dots
\edis
see (\ref{7.2}).
\end{mydef12}

\subsection{}

Next, we give, to be complete, the following

\begin{mydef42}
\be \label{8.4}
\begin{split}
& \prod_{r=1}^{k_2}\left|\frac{\zeta\left(\frac 12+i\alpha_r^{2,2}\right)}{\zeta\left(\frac 12+i\beta_r^2\right)}\right|^2\sim \\
& \sim \frac{\sin^2(\alpha_0^{1,1})}{\cos^2(\alpha_0^{2,2})}
\prod_{r=1}^{k_1}\left|\frac{\zeta\left(\frac 12+i\alpha_r^{1,1}\right)}{\zeta\left(\frac 12+i\beta_r^1\right)}\right|^2+ \\
& + \frac{\cos(2\mu +U)\cos(\mu+U)\cos\mu}{\cos^2(\alpha_0^{2,2})\cos^2(\alpha_0^{3,3})}
\prod_{r=1}^{k_3}\left|\frac{\zeta\left(\frac 12+i\alpha_r^{3,3}\right)}{\zeta\left(\frac 12+i\beta_r^3\right)}\right|^2 ,\ L\to\infty ,
\end{split}
\ee
\be \label{8.5}
\begin{split}
& \prod_{r=1}^{k_1}\left|\frac{\zeta\left(\frac 12+i\alpha_r^{1,1}\right)}{\zeta\left(\frac 12+i\beta_r^1\right)}\right|^2\sim \\
& \sim \frac{\cos^2(\alpha_0^{2,2})}{\sin^2(\alpha_0^{1,1})}
\prod_{r=1}^{k_2}\left|\frac{\zeta\left(\frac 12+i\alpha_r^{2,2}\right)}{\zeta\left(\frac 12+i\beta_r^2\right)}\right|^2- \\
& - \frac{\cos(2\mu +U)\cos(\mu+U)\cos\mu}{\sin^2(\alpha_0^{1,1})\cos^2(\alpha_0^{3,3})}
\prod_{r=1}^{k_3}\left|\frac{\zeta\left(\frac 12+i\alpha_r^{3,3}\right)}{\zeta\left(\frac 12+i\beta_r^3\right)}\right|^2 ,\ L\to\infty .
\end{split}
\ee
\end{mydef42}

\subsection{}

Now, we have the following three sets of oscillating systems
\bdis
\left\{ [\zeta,Q^2,k_1;f_1]\right\}_{k_1=1}^{k_0},\ f_1=\sin^2t,
\edis
\be \label{8.6}
\left\{ [\zeta,Q^2,k_2;f_2]\right\}_{k_2=1}^{k_0},\ f_2=\cos^2t,
\ee
\bdis
\left\{ [\zeta,Q^2,k_3;f_3]\right\}_{k_3=1}^{k_0},\ f_3=\frac{1}{\cos^2t}.
\edis

\begin{remark}
Consequently, new set of interactions corresponds to our three sets of oscillating systems (\ref{8.6}), namely we have the
third order diagram
\bdis
[\zeta,Q^2,k_1;f_1] \longleftrightarrow [\zeta,Q^2,k_2,f_2]\longleftrightarrow [\zeta,Q^2,k_3,f_3]\longleftrightarrow [\zeta,Q^2,k_1,f_1] .
\edis
\end{remark}

\begin{remark}
Now we can expect there are diagrams of the order 4, 5, and so on.
\end{remark}

\section{On new type of factorization formula generated by the second-level elimination}

\subsection{}

Since
\bdis
\int_{2\pi L+\mu}^{2\pi L+\mu+U}\cos t{\rm d}t=2\sin\frac U2\cos\left(\mu+\frac U2\right)
\edis
then
\bdis
\frac 1U\int_{2\pi L+\mu}^{2\pi L+\mu+U}\cos t{\rm d}t=\frac 2U\sin\frac U2\cos\left(\mu+\frac U2\right) ,
\edis
and we obtain by our algorithm the following

\begin{mydef54}
Let
\be \label{9.1}
\begin{split}
& f_4(t)=\cos t,\ t\in [2\pi L+\mu,2\pi L+\mu +U], \\
& U>0,\ \mu\geq\mu_0>0,\ U\in (0,U_0],\ \mu+\frac{U_0}{2}\leq \frac{\pi}{2}-\epsilon,
\end{split}
\ee
where, of course,
\bdis
f_4(t)\in \tilde{C}[2\pi L+\mu,2\pi L+\mu +U].
\edis
Then the following factorization formula holds true
\be \label{9.2}
\prod_{r=1}^{k}\left|\frac{\zeta\left(\frac 12+i\alpha_r^{4}\right)}{\zeta\left(\frac 12+i\beta_r\right)}\right|^2\sim
\frac 2U\sin\frac U2\frac{\cos\left(\mu+\frac U2\right)}{\cos(\alpha_0^4)},\ L\to\infty,
\ee
where
\be \label{9.3}
\begin{split}
& \alpha_0^4=\alpha_0(U,\mu,L,k;f_4), \dots \\
& \beta_r=\beta_r(U,\mu,L,k), \\
& 2\pi L+\mu<\alpha_0^4<2\pi L+\mu +U \ \Rightarrow \ \mu<\alpha_0^4-2\pi L<\mu+U ,
\end{split}
\ee
(comp. (\ref{9.1})).
\end{mydef54}

\subsection{}

Next, we obtain by the similar way the following

\begin{mydef55}
Let
\bdis
\begin{split}
& f_5(t)=\cos t,\ t\in [2\pi L+\mu,2\pi L+\mu +U],
\end{split}
\edis
(under the same assumptions as in (\ref{9.1})).
Then the following factorization formula holds true
\be \label{9.4}
\prod_{r=1}^{k}\left|\frac{\zeta\left(\frac 12+i\alpha_r^{5}\right)}{\zeta\left(\frac 12+i\beta_r\right)}\right|^2\sim
\frac 2U\sin\frac U2\frac{\cos\left(\mu+\frac U2\right)}{\sin(\alpha_0^5)},\ L\to\infty,
\ee
where
\be \label{9.5}
\begin{split}
& \alpha_0^5=\alpha_0(U,\mu,L,k;f_5), \dots \\
& \beta_r=\beta_r(U,\mu,L,k),
\end{split}
\ee
(comp. (\ref{7.3})).
\end{mydef55}

\subsection{}

Since the right-hand side of formulae (\ref{9.2}) and (\ref{9.4}) are bounded and nonzero (see the assumptions in (\ref{9.1})), then we
have the following.

\begin{mydef56}
Under the assumptions (\ref{9.1}) we have the following interaction formula
\be \label{9.6}
\begin{split}
& \prod_{r=1}^{k_2}\left|\frac{\zeta\left(\frac 12+i\alpha_r^{5,2}\right)}{\zeta\left(\frac 12+i\beta_r^2\right)}\right|^2\sim \\
& \sim \tan\left(\mu+\frac U2\right)\frac{\cos(\alpha_0^{4,1})}{\sin(\alpha_0^{5,2})}
\prod_{r=1}^{k_1}\left|\frac{\zeta\left(\frac 12+i\alpha_r^{4,1}\right)}{\zeta\left(\frac 12+i\beta_r^1\right)}\right|^2,\ L\to\infty ,
\end{split}
\ee
where
\be \label{9.7}
\begin{split}
& \alpha_0^{5,2}=\alpha_0^5(U,\mu,L,k_2;f_5),\dots \\
& \beta_r^2=\beta_r(U,\mu,L,k_2), \\
& \alpha_0^{4,1}=\alpha_0^4(U,\mu,L,k_1;f_4),\dots \\
& \beta_r^1=\beta_r(U,\mu,L,k_1), \\
& 1\leq k_1,k_2\leq k_0.
\end{split}
\ee
\end{mydef56}

\subsection{}

Now, in the case
\be \label{9.8}
k_1=k_2=k
\ee
it is true that (see (\ref{9.3}), (\ref{9.7}))
\be \label{9.9}
\beta_r^2=\beta_r^1=\beta_r,\ r=1,\dots ,k_0.
\ee
Next, we have by (\ref{9.9})
\bdis
\prod_{r=1}^{k_2}\left|\zeta\left(\frac 12+i\beta_r^2\right)\right|^2=\prod_{r=1}^{k_1}\left|\zeta\left(\frac 12+i\beta_r^1\right)\right|^2>0
\edis
(comp. also (\ref{1.4}), (\ref{1.5})) and, consequently,
\be \label{9.10}
\begin{split}
&
\frac
{\prod_{r=1}^{k}\left|\frac{\zeta\left(\frac 12+i\alpha_r^{5}\right)}{\zeta\left(\frac 12+i\beta_r\right)}\right|^2}
{\prod_{r=1}^{k}\left|\frac{\zeta\left(\frac 12+i\alpha_r^{4}\right)}{\zeta\left(\frac 12+i\beta_r\right)}\right|^2}=
\prod_{r=1}^{k}\left|\frac{\zeta\left(\frac 12+i\alpha_r^{5}\right)}{\zeta\left(\frac 12+i\alpha_r^4\right)}\right|^2.
\end{split}
\ee

Now, the following theorem holds true.

\begin{mydef13}
Under the assumptions of (\ref{9.1}) we have the following (new) type of interaction formula
\be \label{9.11}
\prod_{r=1}^{k}\left|\frac{\zeta\left(\frac 12+i\alpha_r^{5}\right)}{\zeta\left(\frac 12+i\alpha_r^4\right)}\right|^2\sim
\tan\left(\mu+\frac U2\right)\frac{\cos(\alpha_0^4)}{\sin(\alpha_0^5)},\ L\to\infty,\ k=1,\dots,k_0 .
\ee
\end{mydef13}

\subsection{}

There is a kind of hierarchy in the class of sets of formulae we have obtained:
\begin{itemize}
\item[(a)] First of all, we have two sets (see (\ref{9.2}), (\ref{9.3}))
\be \label{9.12}
\left\{ [\zeta,Q^2,k;f_4]\right\}_{k=1}^{k_0},\ \left\{ [\zeta,Q^2,k;f_5]\right\}_{k=1}^{k_0}
\ee
of the oscillating systems.
\item[(b)] Next, new set containing
\bdis
(k_0)^2
\edis
elements of interactions
\be \label{9.13}
[\zeta,Q^2,k_1;f_4]\longleftrightarrow [\zeta,Q^2,k_2;f_5],\ 1\leq k_1,k_2\leq k_0
\ee
between oscillating systems is generated by the sets (\ref{9.12}) (see formula (\ref{9.6})).
\item[(c)]
Now, by the subset
\be \label{9.14}
[\zeta,Q^2,k;f_4]\longleftrightarrow [\zeta,Q^2,k;f_5],\ k=1,\dots,k_0
\ee
of interactions (\ref{9.13}) is generated the set of new type of factorization formulae (see (\ref{9.11}))
\bdis
\prod_{r=1}^{k}\left|\frac{\zeta\left(\frac 12+i\alpha_r^{5}\right)}{\zeta\left(\frac 12+i\alpha_r^4\right)}\right|^2\sim
\tan\left(\mu+\frac U2\right)\frac{\cos(\alpha_0^4)}{\sin(\alpha_0^5)},\ L\to\infty,\ k=1,\dots,k_0 .
\edis
\end{itemize}

\thanks{I would like to thank Michal Demetrian for his help with electronic version of this paper.}


\begin{thebibliography}{29}
\bibitem{1}
J. Moser,
`Jacob's ladders and almost exact asymptotic representation of the Hardy-Littlewood integral`,
Math. Notes 88, (2010), 414-422, arXiv: 0901.3937.
%
\bibitem{2}
J. Moser,
`Jacob's ladders, structure of the Hardy-Littlewood integral and some new class of nonlinear integral equations`,
Proc. Steklov Inst. 276 (2011), 208-221, arXiv: 1103.0359.
%
\bibitem{3}
J. Moser, `Jacob's ladders, reverse iterations and new infinite set of $L_2$-orthogonal systems generated by the
Riemann zeta-function, arXiv: 1402.2098.
%
\bibitem{4}
J. Moser, `Jacob's ladders, $\zeta$-factorization and infinite set of metamorphoses of a multiform`, arXiv: 1501.07705v2.
%
\bibitem{5}
J. Moser,
`Jacob's ladders, Riemann's oscillators, quotient of the oscillating multiforms and set of metamorphoses of this system`,
arXiv: 1506.00442.
%
\bibitem{6}
J. Moser, `Jacob's ladders, factorization and metamorphoses as an appendix to the Riemann functional equation for $\zeta(s)$ on
the critical line`, arXiv: 1506.00442v1.
%
\bibitem{7}
J. Moser, 'Jacob's ladders, $\mcal{Z}_{\zeta,Q^2}$-transformation of real elementary functions and telegraphic
signals generated by the power functions', arXiv: 1602.04994.
%
\bibitem{8}
C.L. Siegel, `\" Uber Riemann's Nachlass zur analytischen Zahlen-theorie : Quellen und Studien zur Geschichte der Mathematik, Astronomie und Physik`,
Abt. B, Studien 2 (1932), 45--80.



\end{thebibliography}
\end{document}